\documentclass[graybox]{svmult}
\usepackage{listings}
\usepackage{csquotes}    
\usepackage{amsmath}
\usepackage{graphicx}        
\usepackage[bottom]{footmisc}
\usepackage{algorithm,algpseudocode}
\usepackage{pgfplotstable}
\usepackage{pgfplots}
\usepackage{array,multirow,graphicx}
\usepackage{tikz}
\usepgfplotslibrary{groupplots}
\usetikzlibrary{matrix}
\usetikzlibrary{arrows.meta, positioning, chains}
\usetikzlibrary{external}
\usepackage{pgfplots}
\pgfplotsset{compat=newest}
\usepackage{pgfplotstable}
\DeclareMathAlphabet{\pazocal}{OMS}{zplm}{m}{n}    
\pgfplotsset{compat = 1.3}
\usepackage{wrapfig}
\usepackage{xcolor}
\definecolor{TR}{HTML}{6a4c93}
\definecolor{SN}{HTML}{c9182c}

\usepackage{booktabs}

\definecolor{TLTR0}{rgb}{0.63,0.12156,0.48236} 
\definecolor{TLTR1}{HTML}{ff924c}
\definecolor{TLTR2}{HTML}{ffd700}
\definecolor{TLTR3}{HTML}{8ac926}
\definecolor{TLTR4}{HTML}{1982c4}

\definecolor{crimson}{RGB}{220, 20, 60}
\definecolor{darkgreen}{RGB}{0, 150, 40}

\usepackage{color,soul}

\def \0{\mathbf{0}}			%

\def \x{{\boldsymbol{x}}}		
\def \Rm{{\boldsymbol{R}}}	
\def \Bm{{\boldsymbol{B}}}	
\def \Sm{{\boldsymbol{S}}}	
\def \Ym{{\boldsymbol{Y}}}	
\def \Lm{{\boldsymbol{L}}}		
\def \Im{{\boldsymbol{I}}}		
\def \Dm{{\boldsymbol{D}}}	
\def \sv{{\boldsymbol{s}}}	
\def \xv{{\boldsymbol{x}}}	
\def \cv{{\boldsymbol{c}}}	
\def \yv{{\boldsymbol{y}}}		
\def \Cb{{\boldsymbol{C}}}		

\def \R{{\mathbb{R}}}			
\def \Bm{{\mathbf{B}}}			%

\def \l{{\mathbf{l}}}			

\def \I{{\mathbf{I}}}			

\def \x{\vec{x}}				
\def \0{\vec{0}}			%

\def \R{\mathbb{R}}			

\def \I{\mathbf{I}}			

\usepackage{todonotes}
\usepackage{newtxtext}       
\usepackage[varvw]{newtxmath}       

\algnewcommand\algorithmiconput{\textbf{Constants:}}
\algnewcommand\algorithmicinput{\textbf{Input:}}
\algnewcommand\algorithmicoutput{\textbf{Output:}}
\algnewcommand{\algorithmicgoto}{\textbf{go to}}%

\algnewcommand\Constants{\item[\algorithmiconput]}
\algnewcommand\Input{\item[\algorithmicinput]}%
\algnewcommand\Output{\item[\algorithmicoutput]}%
\algnewcommand{\Goto}[1]{\algorithmicgoto~\ref{#1}}%

\AtBeginDocument{%
  \setlength{\abovedisplayskip}{4pt}
  \setlength{\belowdisplayskip}{4pt}
  \setlength{\abovedisplayshortskip}{2pt}
  \setlength{\belowdisplayshortskip}{2pt}
  \setlength{\floatsep}{6pt}        
\setlength{\textfloatsep}{6pt}    
\setlength{\intextsep}{6pt}       
}

\makeindex             

\newcommand{\newdel}[1]{}
\newcommand{\newadd}[1]{#1}

\begin{document}
\title*{Multi-Preconditioned LBFGS for Training Finite-Basis PINNs}
\author{Marc Salvadó-Benasco\orcidID{0009-0001-8647-5159},\\
Aymane Kssim\orcidID{0009-0007-1814-655X},\\
Alexander Heinlein\orcidID{0000-0003-1578-8104},\\
Rolf Krause\orcidID{0000-0001-5408-5271}, \\
Serge Gratton\orcidID{0000-0002-5021-2357}, and\\ 
Alena Kopaničáková\orcidID{0000-0001-8388-5518}.}
\authorrunning{Marc Salvadó-Benasco et al.}

\institute{
Marc Salvadó-Benasco \at Università della Svizzera Italiana, \email{marc.salvado@usi.ch}
\and
 Aymane Kssim \at Toulouse INP-ENSEEIHT, IRIT, ANITI; \email{aymane.kssim@toulouse-inp.fr}
 \and
Alexander Heinlein \at Delft Institute of Applied Mathematics, TU Delft; 
\email{a.heinlein@tudelft.nl} 
\and 
Rolf Krause \at King Abdullah University of Science and Technology;  \email{rolf.krause@kaust.edu.sa} 
\and
 Serge Gratton \at Toulouse INP-ENSEEIHT, IRIT, ANITI; \email{serge.gratton@toulouse-inp.fr}
\and
 Alena Kopaničáková \at Toulouse INP-ENSEEIHT, IRIT, ANITI; \email{alena.kopanicakova@toulouse-inp.fr}
}

\maketitle
\abstract{A multi-preconditioned LBFGS (MP-LBFGS) algorithm is introduced for training finite-basis physics-informed neural networks (FBPINNs).
The algorithm is motivated by the nonlinear additive Schwarz method and exploits the domain-decomposition-inspired additive architecture of FBPINNs, in which local neural networks are defined on subdomains, thereby localizing the network representation. Parallel, subdomain-local quasi-Newton corrections are then constructed on the corresponding local parts of the architecture.
A key feature is a novel nonlinear multi-preconditioning mechanism, in which subdomain corrections are optimally combined through the solution of a low-dimensional subspace minimization problem.
Numerical experiments indicate that MP-LBFGS can improve convergence speed, as well as model accuracy over standard LBFGS while incurring lower communication overhead.

}


\newcommand{\NN}{\pazocal{N}}
\newcommand{\NNj}{\pazocal{N}_j}
\renewcommand{\l}{\pazocal{L}}
\newcommand{\lphy}{\l_{\text{phy}}}
\newcommand{\lbc}{\l_{\text{bc}}}
\newcommand{\ldata}{\l_{\text{data}}}
\newcommand{\lamphy}{\lambda_{\text{phy}}}
\newcommand{\lambc}{\lambda_{\text{bc}}}
\newcommand{\lamdata}{\lambda_{\text{data}}}

\newcommand{\ns}{n_s}
\newcommand{\nb}{n_b}
\newcommand{\xii}{x^{(i)}}
\newcommand{\xiib}{\x^{(i)}}
\newcommand{\X}{\mathbf{X}}
\newcommand{\Xj}{\X_j}

\renewcommand{\D}{\pazocal{D}}
\newcommand{\A}{\pazocal{A}}
\newcommand{\V}{\pazocal{V}}
\newcommand{\Vj}{\V_j}

\newcommand{\norm}[1]{\left\lVert\, #1 \, \right\rVert}
\newcommand{\inorm}[1]{\norm{#1}_{\infty}}
\newcommand{\diag}{\operatorname{diag}}

\newcommand{\oned}[1]{\mathbf{1}_{#1}}

\newcommand{\thetab}{\boldsymbol{\theta}}
\section{Introduction}

Many scientific and engineering applications require the solution of partial differential equations (PDEs).
Classical numerical methods, such as finite element (FE) discretizations, are accurate but can become computationally prohibitive for high-dimensional, multiscale, or data-driven problems. Neural network-based discretization methods, such as
Physics-informed neural networks (PINNs)~\cite{lagaris1998artificial, raissi2019physics}, offer a mesh-free alternative that can naturally incorporate observational data by enforcing physical constraints through the training loss. This comes at the cost of a challenging, non-convex optimization problem for training the PINN. While standard PINNs struggle with multiscale or highly oscillatory solutions, finite-basis PINNs (FBPINNs)~\cite{moseley2023finite,dolean2024multilevel} mitigate these difficulties through an additive, domain-decomposition (DD) inspired architecture defined on overlapping subdomains, in which collocation points are split and shared across neighboring subdomains.

The training of PINNs, including FBPINNs, has predominantly been carried out using stochastic gradient descent (SGD) and its variants. Alternatively, quasi-Newton methods, such as the limited-memory Broyden-Fletcher-Goldfarb-Shanno (LBFGS) algorithm~\cite{liu1989limited}, have been considered due to their ability to incorporate curvature information~\cite{kiyani2025optimizer}.
In FBPINNs, the additive structure of the model allows the forward and backward passes to be computed in parallel across subdomains. However, a parallel implementation of the LBFGS optimizer still requires global synchronization at each iteration, which can lead to communication-dominated execution and limit overall parallel scalability.
In this work, we propose to address this difficulty by introducing a multi-preconditioned LBFGS (MP-LBFGS) method that enables multiple local optimization steps prior to global synchronization, thereby increasing local work while reducing the communication overhead.

The proposed MP-LBFGS algorithm is motivated by the nonlinear additive Schwarz method~\cite{cai2002nonlinearly} and explicitly exploits the FBPINN decomposition of both the computational domain and the model parameters.
The method performs a sequence of local LBFGS iterations on each subdomain, followed by a global LBFGS synchronization step.
To aggregate the resulting local corrections, we introduce a subspace minimization strategy for scaling and combining subdomain updates.
This is particularly essential in the FBPINN setting, where standard scaling strategies from classical DD may be ineffective due to fundamental structural differences between neural network models and FE discretizations.
As demonstrated by our numerical experiments, this novel scaling mechanism enables stable aggregation of local updates while preserving the benefits of localized desynchronized optimization.

This manuscript is organized as follows.
Sect.~\ref{sec:fbpinn} introduces the FBPINN architecture, Sect.~\ref{sec:MPLBFGS} presents the multi-preconditioned LBFGS algorithm, and numerical results and conclusions are given in Sect.~\ref{sec:num_res} and Sect.~\ref{sec:conclusion}, respectively.

\section{Finite Basis PINNs (FBPINNs)} 
\label{sec:fbpinn}
PINNs are neural network (NN) models that approximate solutions of differential equations.
To this aim, let $\Omega \subset \mathbb{R}^d$ be a bounded domain with boundary $\partial \Omega$.
We then consider the following abstract boundary value problem:
\begin{equation}
\begin{aligned}
\pazocal{P}[u](\x) &= f(\x),
&\qquad &\forall \x \in \Omega, \\
u(\x) &= g(\x),
&\qquad & \forall \x \in \partial\Omega, 
\label{eq:model_problem}
\end{aligned}
\end{equation}
where $\pazocal{P}$ denotes a differential operator,
and $f$ is a forcing term. Note that other boundary conditions can be treated analogously.

We aim to approximate the unknown solution $u : \Omega \to \mathbb{R}$ using an NN model
$\NN: \mathbb{R}^p  \times \Omega \to \mathbb{R}$,
which is parameterized by a set of weights $\thetab \in \mathbb{R}^p$, where $p$ denotes the number of trainable parameters.
To determine suitable network parameters, we sample a set of $n$ collocation points $\pazocal{D} := \{ \x_i \}_{i=1}^{n} \subset \Omega$ located in the interior of $\Omega$.
In addition, we consider a boundary dataset
$\pazocal{D}_{\mathrm{bc}}$, with $n^{BC}$ collocation points sampled on the boundary $\partial\Omega$. The optimal network parameters are obtained by solving the following nonlinear minimization problem:
\begin{align}
\label{eq:thetastar}
\thetab^* :=
\arg\min_{\thetab} \;
\lamphy \, \lphy(\thetab; \pazocal{D})
+ \lambc \, \lbc(\thetab; \pazocal{D}_{\mathrm{bc}}),
\end{align}
where $\lamphy, 
\lambc \in \mathbb{R}^+$ 
are weighting parameters. The physics loss
$\lphy : \mathbb{R}^{p} \times \pazocal{D} \to \mathbb{R}$
quantifies the violation of the governing differential equation, i.e., 
\[
\lphy(\thetab; \pazocal{D}) :=
\frac{1}{n} \sum_{i=1}^{n}
\Big(
\pazocal{P}\big[\hat{u}(\thetab; \x_i)\big] - f(\x_i)
\Big)^2,
\]
where $\hat{u}$ denotes the NN approximation of the solution.
Similarly, each boundary loss
$\lbc : \mathbb{R}^{p} \times \pazocal{D}_{bc} \to \mathbb{R}$
enforces the boundary conditions and is defined as
\[
\lbc(\thetab; \pazocal{D}_{bc}) :=
\frac{1}{n^{BC}}
\sum_{j}^{n^{BC}}
\Big( \hat{u}(\thetab; \x_j) - g(\x_j)
\Big)^2.
\]

The weights $\lamphy$ and $\lambc$ in~\eqref{eq:thetastar} are used to balance the contributions of the different loss terms.
In practice, however, it can be challenging to choose these parameters so that the boundary conditions (BCs) are enforced with sufficient accuracy while simultaneously ensuring a satisfactory decrease of the physics loss~\cite{lagaris1998artificial}.
To overcome this difficulty, we enforce BCs by designing the solution ansatz such that they are automatically satisfied~\cite{lagaris1998artificial}.
By choosing a function $\ell: \Omega \to \mathbb{R}$ which vanishes on the boundary $\partial \Omega$, we define the network solution approximation as 
\begin{equation} \label{eq:lifting}
    \hat{u}(\thetab; \x)
    :=
    C(\x) + \ell(\x)\, \NN(\thetab; \x),
\end{equation}
where the function $C(\x)$ encodes the prescribed BC through data $g$. With this construction, the BC are satisfied by design, 

and~\eqref{eq:thetastar} can  be reformulated as
\begin{equation} \label{eq:optprob}
\thetab^* :=
\arg\min_{\thetab} \, \l(\thetab; \pazocal{D})
:= \lphy(\thetab; \pazocal{D}).
\end{equation}


\subsection{FBPINN architecture}
Similarly to standard neural networks, PINNs also suffer from spectral bias~\cite{moseley2023finite}, meaning that they tend to learn low-frequency components efficiently while struggling to capture high-frequency features.
To address this limitation, a domain-decomposition-motivated architecture, termed the finite-basis PINN (FBPINN), was introduced in~\cite{moseley2023finite}.
FBPINNs decompose the computational domain into several overlapping subdomains, within which high-frequency components are effectively rescaled to lower frequencies, thereby improving their resolution during training.

Let us decompose the domain $\Omega$ into $\ns$ subdomains, such that $\bigcup_{j=1}^{\ns} \Omega_j = \Omega$, and the width of the overlap between neighboring subdomains is denoted by the parameter $\delta > 0$; see also~\cite{dolean2015introduction,moseley2023finite}. 

For each subdomain $\Omega_j$, we introduce a dataset $\pazocal{D}_j$ of collocation points, with
$\pazocal{D} = \bigcup_{j=1}^{\ns} \pazocal{D}_j$.

\begin{figure}[t]
  \centering
  \includegraphics{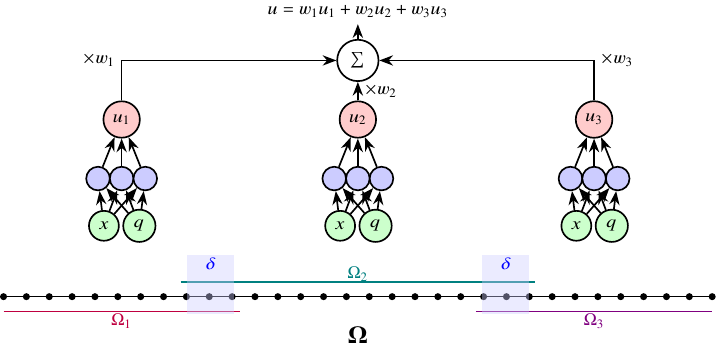}

  \caption{FBPINN with collocation points decomposed into three overlapping subdomains.
    Each subdomain is associated with a subnetwork, and each subnetwork takes as input the spatial coordinates $x$ and, possibly, additional physical features $q$.}
  \label{fig:tikz-external}
\end{figure}

For each subdomain $\Omega_j$, we define a corresponding space of functions as
\begin{equation} \label{eq:vj}
\Vj := \{ \NNj \mid \NNj : \mathbb{R}^{p_j} \times \Omega  \to \mathbb{R} \},
\end{equation}
where $\NNj$ denotes a DNN associated with the $j$-th subdomain and parameterized by $p_j$ trainable parameters.
To restrict the support of each local network $\NNj$ to its corresponding subdomain $\Omega_j$ and to appropriately weight the overlap between subdomains, we introduce a collection of ``window'' functions
$\{ w_j \}_{j=1}^{\ns} : \Omega \to \mathbb{R}$ such that
$\operatorname{supp}(w_j) \subset \Omega_j$ for all $j \in \{1,\dots,\ns\}$ and
$\sum_{j=1}^{\ns} w_j \equiv 1$ on  $\Omega$; this means that they form a partition of unity. \newdel{, see also Fig.~1.}
Then, we define the global approximation space as $\V := \sum_{j=1}^{\ns} w_j \, \Vj$.
The solution $u$ is then approximated by a global network $\NN$, obtained by combining the predictions of all subnetworks, i.e., 
\begin{equation} \label{eq:fbpinn}
\NN(\thetab; \x)
=
\sum_{j=1}^{\ns}
w_j(\x)\,
\NNj\big(\thetab_j; \operatorname{norm}_j(\x) \big).
\end{equation}
Here, the global parameter vector is given by $\thetab = (\thetab_1, \ldots, \thetab_{\ns})$. 
We use restriction operator $\Rm_j:\R^p \to \R^{p_j}$ to extract $\thetab_j$ from $\thetab$, i.e., $\thetab_j = \Rm_j \thetab$.
Conversely, $\Rm_j^\top$ is used to assign local parameters to the global network, i.e., $\thetab = \sum_{j} \Rm_j^\top \thetab_j$. 
Note, there is no overlap between the parameters of different subnetworks, see also Fig.~\ref{fig:tikz-external}. 
Furthermore, to mitigate the spectral bias, the normalization function $\operatorname{norm}_j : \Omega \to (-1,1)^d$ is used to rescale the input coordinates on $\Omega_j$ in order to map high-frequency features in the global domain to lower-frequency representations on each subdomain.
The FBPINN is then trained by inserting~\eqref{eq:fbpinn} into the loss function~\eqref{eq:optprob}; cf.~\cite{moseley2023finite,dolean2024multilevel}.

\section{Multi-preconditioned LBFGS (MP-LBFGS)}
\label{sec:MPLBFGS}
Quasi-Newton (QN) methods are widely used for training PINNs, as they exploit curvature information without requiring the explicit computation or storage of the Hessian matrix.
At each iteration $k$, the network parameters are updated as
\begin{equation}
\thetab^{(k+1)} = \thetab^{(k)} - \alpha^{(k)} (\Bm^{(k)})^{-1} \nabla \pazocal{L}(\thetab^{(k)}),
\label{eq:QN_update}
\end{equation}
where $\Bm^{(k)}$ is a low-rank Hessian approximation and $\alpha^{(k)}$ is obtained by line-search method with strong Wolfe's conditions~\cite{wolfe1969convergence}. 

Traditionally, approximation $\Bm^{(k)}$ is constructed to satisfy the secant equation:
\[
\Bm^{(k+1)} \, \sv^{(k)} = \yv^{(k)},
\]
where $\sv^{(k)} := \thetab^{(k+1)} - \thetab^{(k)}$ and
$\yv^{(k)} := \nabla_{\thetab} \l(\thetab^{(k+1)}) - \nabla_{\thetab} \l(\thetab^{(k)})$.
Since the secant equation alone does not uniquely determine $\Bm^{(k+1)}$, additional conditions must be imposed, leading to different variants of QN methods.
In this work, we focus on the Broyden-Fletcher-Goldfarb-Shanno (BFGS) method~\cite{nocedal1980updating}, as it is the most widely used in the context of PINNs~\cite{kiyani2025optimizer}.

In particular, we consider the limited-memory variant (LBFGS)~\cite{nocedal1980updating}, which approximates the Hessian using only the last $Q$ secant pairs ${ (\sv_{k-q}, \yv_{k-q}) }_{q=1}^{Q}$.
A standard compact representation of the LBFGS is given as
\[
\Bm^{(k+1)} = \Bm^{(0)} - \begin{bmatrix}
\Bm^{(0)}\Sm^{(k)} & \Ym^{(k)}
\end{bmatrix}
\begin{bmatrix}
(\Sm^{(k)})^{\top}\Bm^{(0)}\Sm^{(k)} & \Lm^{(k)} \\
(\Lm^{(k)})^{\top} & -\Dm^{(k)}
\end{bmatrix}^{-1}
\begin{bmatrix}
(\Sm^{(k)})^{\top}\Bm^{(0)}\\
(\Ym^{(k)})^{\top}
\end{bmatrix},
\]
where $\Sm^{(k)} := [\, \sv^{(k-Q+1)}, \ldots, \sv^{(k)} \,]$ and $\Ym^{(k)} := [\, \yv^{(k-Q+1)}, \ldots, \yv^{(k)} \,]$, for ${k>1}$.
Note, that the secant pair $(\sv^{(k)}, \yv^{(k)})$ is added to the L-BFGS memory only if $(\yv^{(k)})^\top \sv^{(k)} > 0$.
Otherwise, the pair is discarded to preserve positive definiteness of the Hessian approximation~\cite{nocedal1980updating}.
The symbols $\Lm^{(k)}$ and $\Dm^{(k)}$ denote the strictly lower triangular and diagonal parts of $(\Ym^{(k)})^{\top} \Sm^{(k)}$, while $\Bm^{(0)}$ is an initial  approximation, e.g.,~$\Bm^{(0)} = \gamma \Im$, where $\gamma >0$. 

Note, to compute the search direction in~\eqref{eq:QN_update}, $(\Bm^{(k)})^{-1}$ is required.
Rather than solving the associated linear system, the product $(\Bm^{(k)})^{-1} \nabla \pazocal{L}(\thetab^{(k)})$ can be computed efficiently using, for example, the well-known two-loop recursion algorithm~\cite{nocedal1980updating}. 

\subsection{Nonlinearly preconditioned LBFGS}
\label{sec:search_dir}
Minimizing the loss $\pazocal{L}$ leads to the following first-order optimality condition:
\begin{equation}
\nabla \l(\thetab; \pazocal{D}) = \mathbf{0}.
\label{eq:kkt}
\end{equation}
To accelerate convergence of LBFGS, we introduce a nonlinear preconditioning operator
$P : \mathbb{R}^p \to \mathbb{R}^p$ that approximates the inverse of the gradient mapping, i.e.,
$P \approx \nabla \l^{-1}$.

While such an operator cannot be constructed explicitly, it can be defined implicitly via a nonlinear fixed-point map of the form $\thetab = P(\thetab; \pazocal{D})$.

Following a right-preconditioning strategy~\cite{Kopanicakova2024, kothari2022nonlinear}, we reformulate the optimality system by introducing the preconditioned gradient as
$\pazocal{F}(\thetab; \pazocal{D}):=\nabla \l\big(P(\thetab; \pazocal{D}); \pazocal{D}\big).$
The original nonlinear problem~\eqref{eq:kkt} is thus replaced by $\pazocal{F}(\thetab; \pazocal{D}) = \mathbf{0}$, which we solve using an LBFGS method.
Note, the preconditioner $P$ induces a nonlinear change of coordinates in the parameter space, while the quasi-Newton method operates on $\pazocal{F}$.

At the $k$-th iteration, the right-preconditioned LBFGS method proceeds in two steps.
First, we perform a nonlinear change of variables, thus the current iterate $\thetab^{(k)}$ is mapped through the nonlinear preconditioner, which defines an intermediate (preconditioned) iterate, i.e., ${\thetab^{(k+1/2)} = P(\thetab^{(k)}; \pazocal{D})}$. 
Second, a quasi-Newton correction is applied to the preconditioned gradient, yielding
\begin{equation}
\thetab^{(k+1)}
=
\thetab^{(k+1/2)}
-
\alpha^{(k)}
\big(\Bm_{\pazocal{F}}^{(k)}\big)^{-1}
\pazocal{F}\big(\thetab^{(k+1/2)}; \pazocal{D}\big).
\end{equation}
Here, the update is taken around the preconditioned iterate $\thetab^{(k+1/2)}$ and the LBFGS correction is applied in the preconditioned coordinates associated with $\pazocal{F}$.
The approximate Hessian $\Bm_{\pazocal{F}}^{(k)}$ is constructed using the standard LBFGS update with secant pairs $\sv^{(k)} = \thetab^{(k+1)} - \thetab^{(k+1/2)}$, and  $\yv^{(k)} =
\pazocal{F}(\thetab^{(k+1)}) -
\pazocal{F}(\thetab^{(k+1/2)})$.

To define the fixed-point operator $P$,  we exploit the FBPINN structure, i.e., 
\begin{equation}
P(\thetab^{(k)}; \pazocal{D})
=
\thetab^{(k)}
+
\beta^{(k)}
\sum_{j=1}^{n_s}
\Rm_j^{\top}
\big(
\thetab_j^{(\ast)} - \Rm_j \thetab^{(k)}
\big),
\label{eq:prec_step}
\end{equation}
where $\beta^{(k)} \in \R^+$ and
$\thetab_j^{(\ast)}$ is a solution of the local minimization problem
\begin{equation}
\thetab_j^{(\ast)}
=
\arg\min_{\thetab_j}
\pazocal{L}(\thetab; \pazocal{D}),
\label{eq:subproblem_min}
\end{equation}
where $\thetab_j$ denotes the set of trainable parameters of the j-th subdomain.

Eq.~\eqref{eq:prec_step} 
corresponds to a nonlinear additive Schwarz preconditioner in parameter space.
In practice, the local problem~\eqref{eq:subproblem_min} does not need to be solved exactly.
Instead, an approximate solution suffices and can be obtained, for instance, by performing a fixed number $\eta$ of L-BFGS iterations initialized at $\thetab_j^{(k)} = \Rm_j \thetab^{(k)}$.
Here, we remark that local and global LBFGS maintain separate secant memories. 
The resulting right-preconditioned LBFGS algorithm is summarized in Alg.~\ref{alg:multi_precon_lbfgs}. 

\begin{algorithm}[t]
\caption{Multi-Preconditioned LBFGS (MP-LBFGS)}
\label{alg:multi_precon_lbfgs}
\begin{algorithmic}[1]
\Require $\boldsymbol{\theta}^{(0)} \in \R^p$, $\pazocal{D}$, $\pazocal{N}: \R^p \times \Omega \rightarrow \R$, $k_{\max} \in \mathbb{N}$

\For{$k=0, \ldots, k_{\max}$}
    \For{$j \in \{1,...,n_{s}\}$} \Comment{Parallel preconditioning step}

        \State $\boldsymbol{\theta}^{(k + 1/2)}_j \gets \text{LBFGS}(\pazocal{N}_j,  \mathbf{R}_j \boldsymbol{\theta}^{(k)}, \eta)$ \Comment{Perform $\eta$ steps of local LBFGS}
        \State $\mathbf{c}_j \gets \boldsymbol{\theta}^{(k + 1/2)}_j -  \mathbf{R}_j \boldsymbol{\theta}^{(k)}$
    \EndFor

    \State $\boldsymbol{\theta}^{(k + 1/2)} \gets \boldsymbol{\theta}^{(k)} + \sum_{j=1}^{n_{s}} {\beta}^{(k)}_j \mathbf{R}_j^\intercal \mathbf{c}_j$  \Comment{Compute $\boldsymbol{\beta}^{(k)}$ using approaches from Sect.~\ref{sec:multi_prec} 
    }
    \State $\boldsymbol{\theta}^{(k + 1)}= \text{LBFGS}(\pazocal{N}, \boldsymbol{\theta}^{(k + 1/2)}, 1)$ \Comment{Perform one step of global LBFGS}
\EndFor
\end{algorithmic}
\end{algorithm}

\subsection{Scaling of Subdomain Corrections} 
\label{sec:multi_prec}
The preconditioning step~\eqref{eq:prec_step} yields a set of directions $\{\mathbf{c}_j\}_{j=1}^{n_{s}}$, where each $\mathbf{c}_j := \thetab_j^{(\ast)} - \Rm_j \thetab^{(k)}$.
In this section, we discuss how to optimally update the global parameters using these search directions. 
Therefore, at each iteration $k$, we are looking for scaling parameters
$\boldsymbol{\beta} = \begin{bmatrix}
            \beta_1,
            \beta_2,
            \hdots ,
            \beta_{n_{s}}
        \end{bmatrix}^\intercal,$
such that
\begin{equation} \label{eq:update_condition}    
    \pazocal{L}\bigg(\mathbf{\boldsymbol{\theta}} +\sum_{j = 1}^{n_{s}} \beta_j \mathbf{R}_j^{\intercal}\mathbf{c}_j\bigg) <
    \pazocal{L}\big(\mathbf{\boldsymbol{\theta}}\big).
\end{equation}
We emphasize that the choice of the scaling parameters $\boldsymbol{\beta}$ is critical due to the intrinsic structure of FBPINNs and their differences from classical DD methods.
\newadd{In particular, classical DD methods, originally proposed for elliptic PDEs, exploit structural properties of finite-element/difference discretizations, such as locality of basis functions and sparse, spatially localized coupling between degrees of freedom. 
This ensures that subdomain corrections mainly affect local error components, while global modes are handled through coarse-grid or global synchronization steps.}

\newadd{In contrast, machine-learning problems involve nonconvex minimization, for which only a few domain-decomposition methods have been developed, e.g.,~\cite{gratton2025recursive}. 
Moreover, locally defined physical degrees of freedom are replaced by globally coupled parameters, fundamentally altering error-propagation mechanisms.
As a result, classical DD assumptions, such as convexity, locality, and spectral separation, no longer hold, which often leads to degraded convergence of conventional DD methods.}

\newdel{Although the subnetworks do not share parameters (zero overlap), overlap in the spatial domain induces indirect coupling through their outputs, which causes the desynchronization of the solutions during the preconditioning step.
In this work, we investigate three strategies for selecting the scaling parameters.}
\newadd{For the FBPINNs considered in this work, spatial overlap induces indirect coupling, leading to desynchronization of parameters during the preconditioning step and motivating the design of novel scaling strategies for the subdomain corrections. 
In particular, we investigate three such strategies:}

\begin{table}[t]

    \centering
    {
        \caption{Estimates of the parallel number of loss function evaluations ($\# \pazocal{L}_{\text{e}}$), gradient evaluations ($\# g_{\text{e}}$), update cost (UC), and memory cost (MC) per iteration.}
        \label{tab:computational_cost}
        \begin{tabular}{|l|l|l|l|l|l|}
            \hline
            Method                           & $\# \pazocal{L}_{\text{e}}$                                   & $\# g_{\text{e}}$              & UC                                                                          & MC                              \\ [0.5ex]
            \hline
            \hline
            LBFGS                           & $1 + \text{its}_{ls}$                                         & $1$                            & $p+ 4 Q p$                                                                  & $p+ Q p$                        \\
            \hline
           MP-LBFGS   (UniS)  & $1 + \text{its}_{ls} +    \eta (\# \pazocal{L}_{\text{e}_{j}})$                                      & $2 + \eta (\# g_{\text{e}_{j}})$                         & $2p+4 Q p+ \eta \text{UC}_{j}$                                                                 & $2p+2 Q p+ \text{MC}_{j}$                     \\
            \hline
             & $1 + {\#\text{its}_{\text{ls}}} + $                           & $2 + $                         & $2p+4 Q p+ $                                                                & $2p+2 Q p+$                     \\
               MP-LBFGS                               & $  \eta \# \pazocal{L}_{\text{e}_{j}} + 2 + $                     & $ \eta \# g_{\text{e}_{j}} + 2 +$ & $\eta \text{UC}_{j} + 2\frac{p}{n_{s}} +$                                        & $\text{MC}_{j}$ +   $2p +$ \\
                            (SPM)                       & $\#\text{its}_{\text{Newton}} (1 + \#\text{its}_{\text{ls}})$ & $\#\text{its}_{\text{Newton}}$ & $\#\text{its}_{\text{Newton}} (\frac{p}{n_{s}} + \#\text{its}_{\text{ls}})$ & $2n_{s} + n_{s}^2$              \\
            \hline
        \end{tabular}
    }
\end{table}

\begin{itemize}
    \item \textbf{Uniform scaling (UniS):} 
    As a first approach, we consider uniform scaling by a constant $\beta_0$, i.e.,
$\beta_j = \beta_0$, for all $j = 1, \ldots, n_s$.
The value of $\beta_0$ is typically set to $1$ in non-overlapping DD methods, but it can be fine-tuned.
This approach has a low computational cost and is naturally parallel.
However, it has been shown in~\cite{lee2026two} that it does not yield an efficient preconditioner for NN problems.

    \item \textbf{Multiplicative line search scaling (LSS):} 
Following~\cite{lee2026two}, we also investigate a line search approach,
    where we start with a search direction $\mathbf{c}_1$ on the first subdomain, and look for the biggest value $\beta_1$ that satisfies  $\pazocal{L} (\mathbf{\boldsymbol{\theta}} + \beta_1 \mathbf{R}_1^{\intercal}\mathbf{c}_1) < \pazocal{L} (\mathbf{\boldsymbol{\theta}})$.
    We then iterate through all search directions (subdomains) and seek their corresponding scaling factors such that
    \begin{equation}
        \pazocal{L}\bigg(\mathbf{\boldsymbol{\theta}} + \sum_{j=1}^{i-1}\beta_j \mathbf{R}_j^{\intercal}\mathbf{c}_j +\beta_i \mathbf{R}_i^{\intercal}\mathbf{c}_i \bigg) <
        \pazocal{L}\bigg(\mathbf{\boldsymbol{\theta}} + \sum_{j=1}^{i-1}\beta_j \mathbf{R}_j^{\intercal}\mathbf{c}_j\bigg),
        \quad \forall i = {2,..., n_{s}}.
    \end{equation}
In practice, to ensure a sufficient decrease of the loss, we employ a line search~\cite{wolfe1969convergence}.

The LSS scheme introduces two challenges: first, its inherently sequential, and second, the scheme implicitly depends on the ordering of the search directions.

    \item \textbf{Subspace minimization (SPM):} 
Given the limitations of the previous approaches, we propose determining the subdomain scaling parameters $\boldsymbol{\beta}$ by solving a low-dimensional subspace minimization problem.
Motivated by linear multi-preconditioning strategies~\cite{bridson2006multipreconditioned}, we define
$ \Cb :=
\begin{bmatrix}
\mathbf R_1^{\intercal} \cv_1 &
\mathbf R_2^{\intercal} \cv_2 &
\cdots &
\mathbf R_{n_s}^{\intercal} \cv_{n_s}
\end{bmatrix}$.

The scaling vector $\boldsymbol{\beta} \in \mathbb{R}^{n_s}$ can now be obtained by solving the subspace problem
\begin{equation} \label{eq:subproblem}
\min_{\boldsymbol{\beta}} \;
\phi(\boldsymbol{\beta})
:=
\pazocal{L}(\boldsymbol{\theta} + \Cb \boldsymbol{\beta}).
\end{equation}

To solve~\eqref{eq:subproblem}, we perform a few iterations of a simplified Newton method~\cite{deuflhard2011newton}, updating $\boldsymbol{\beta}$ at iteration $m$ using the following update rule:
\[
\boldsymbol{\beta}^{(m+1)}
=
\boldsymbol{\beta}^{(m)}
-
t^{(m)}
\big(\nabla^2 \phi(\boldsymbol{\beta}^{(0)})\big)^{-1}
\nabla \phi(\boldsymbol{\beta}^{(m)}),
\]
where $t^{(m)}$ is a damping parameter chosen by line search~\cite{wolfe1969convergence}.
The gradient and Hessian of $\phi$ are given by
$
\nabla \phi(\boldsymbol{\beta})
=
\Cb^{\intercal} \nabla \pazocal{L}(\boldsymbol{\theta} + \Cb \boldsymbol{\beta})$, and
$
\nabla^2 \phi(\boldsymbol{\beta})
=
\Cb^{\intercal} \nabla^2 \pazocal{L}(\boldsymbol{\theta} + \Cb \boldsymbol{\beta}) \Cb,
$
respectively. 
To avoid explicit Hessian computations, we approximate the action of
$\nabla^2 \pazocal{L}(\boldsymbol{\theta})$ on $\Cb$ using finite differences and precompute
$\big(\nabla^2 \phi(\boldsymbol{\beta}^{(0)})\big)^{-1}$ at $m=0$.
Since $\nabla^2 \phi$ has size $n_s \times n_s$, the resulting system remains inexpensive to solve.
Moreover, although each Newton step requires multiple forward and backward evaluations, these operations are fully parallel across subdomains.    

\end{itemize}

Tab.~\ref{tab:computational_cost} 
compares the estimated parallel computational costs of LBFGS and MP-LBFGS.
Here, the word parallel refers to the per-device computational cost assuming a parallel implementation of FBPINN (one subdomain, one subnetwork per device), LBFGS, and MP-LBFGS.
Notably, communication cost is not included in this estimate.
For LBFGS, each epoch involves one forward and one backward pass to evaluate the gradient.
In addition, the line search in~\eqref{eq:QN_update} requires $\#\text{its}_{\text{ls}}$ forward evaluations.
The computational cost of the LBFGS update scales with the number of stored secant pairs $Q$ as $2p + 4Qp$, while the memory requirements accounting for the parameters, secant pairs, and momentum terms scale as $2p + 2Qp$.

The cost of MP-LBFGS depends on the choice of scaling strategy.
In all cases, $\eta$ forward and backward passes are required for each subnetwork, denoted by $\# \pazocal{L}{\text{e}_{j}}$ and $\# g_{\text{e}_{j}}$, respectively.
Since the training of each subnetwork can be performed in parallel, the local L-BFGS update cost is denoted by $\text{UC}_{j}$, and the corresponding memory requirements by $\text{MC}_{j}$.
The cost of global L-BFGS step corresponds to that of LBFGS, with an addition of computing and storing the gradient at $\thetab^{(k+1/2)}$.

For SPM, two additional loss and gradient evaluations are needed to compute $\nabla^2 \pazocal{L} (\boldsymbol{\theta}) \Cb$. 
Moreover, at each Newton iteration, we require gradient evaluation to compute $\Cb^{\intercal} \nabla \pazocal{L}(\boldsymbol{\theta} + \Cb \boldsymbol{\beta})$ as well as $\#\text{its}_{\text{ls}}$ loss calls to  find $t^{(m)}$.
Storing the quantity $\nabla^2 \pazocal{L}(\boldsymbol{\theta}) \Cb$ requires storing $p$ parameters per compute node.
Each product $\nabla^2 \pazocal{L}(\boldsymbol{\theta}) \mathbf R_j^{\intercal} \mathbf c_j$ can be computed and stored locally, after which $\nabla^2 \phi$ is assembled by multiplying with $\Cb^{\intercal}$ and concatenating the resulting columns.
In addition, the vectors $\boldsymbol{\beta}^{(m)}$ and Newton's search direction, as well as the Hessian matrix $\nabla^2 \phi \in \mathbb{R}^{n_s \times n_s}$, must be stored.
The UC further includes cost associated with Hessian evaluation.

\section{Numerical results} \label{sec:num_res}
\begin{figure}[t]
    \centering
    \includegraphics{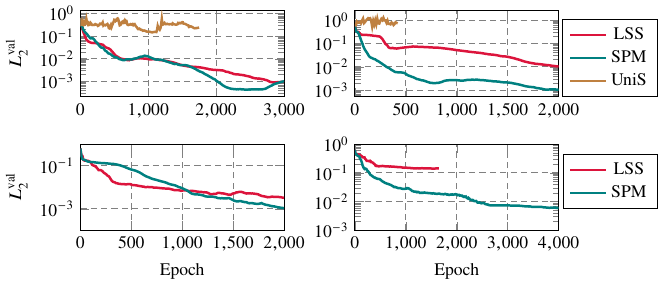}

    \caption{Convergence of relative error $L_2^{\text{val}}$ for Burger's equation. 
    Problem is trained using MP-LBFGS with different subdomains scaling strategies. 
    We consider $4 \times 2$ (top left), $4 \times 4$ (top right), $2 \times 2$ (bottom left), and $4 \times 1$ (bottom right) subdomains.}
    \label{fig:plotscalings_robustness}
\end{figure}

In this section, we investigate the numerical performance of our MP-LBFGS algorithm.
To this aim, we consider two benchmark problems, namely:
\begin{itemize}
\item \textbf{Poisson's equation:} 

We consider one- and two-dimensional Poisson equations with homogeneous Dirichlet boundary conditions, i.e., 
\begin{equation}
	\begin{aligned}
		-\Delta u &= f, & \quad &\forall \xv \in \Omega, \\
		u &= 0, & \quad &\forall \xv \in \partial \Omega.
	\end{aligned}
	\label{eq:poisson}
\end{equation}
For $\Omega = (0,1)^2$, $f$ is chosen such that the exact solution is
$u_{\mathrm{true}}(\xv) = \sin(4\pi \xv_1)\sin(4\pi \xv_2)$.
For $\Omega = (0,1)$, the exact solution is $u_{\mathrm{true}}(\xv) = \sin(20\pi \xv_1)$.
Validation is performed using the relative $L_2$ error with respect to $u_{\mathrm{true}}$.

\item \textbf{Burgers' equation:}
Let $\Omega = (0, 1) \times (-1, 1)$, we consider the Burgers' equation:
\begin{equation}
	\begin{aligned}
		\frac{\partial u}{\partial t} + u \nabla u - \nu \nabla^2 u & = 0,  \quad \quad \quad &  & \forall \ (t,x) \in (0,1] \times (-1,1), \\
		u(0, x)                  & = - \sin(\pi x)   \                &  & \forall x \in [-1,1],                    \\
		u(t, 1) = u(t, -1)    & = 0 \                   &  & \forall \ t \in (0,1],                    \\
	\end{aligned}
\end{equation}
where $\nu = 0.01 / \pi$ is the kinematic viscosity.
The accuracy of the FBPINN solution is assessed by comparison with a finite element reference solution obtained on a mesh with $25{,}600$ degrees of freedom, using the relative $L_2$ error as the metric.
\end{itemize}

All problems are posed on uniform domains, with 2D decompositions performed uniformly in both directions.
The weighting functions are 
$
	w_j(\mathbf{x}) := \prod_{k=1}^d \Big(1 + \cos\big(\pi\, \operatorname{norm}_j(\mathbf{x})|_k\big)\Big)^2.
$
Each subnetwork is a four-layer ResNet~\cite{he2016deep} of width $20$ with $\tanh$ activation.
The datasets consist of $20{,}000$ collocation points for 2D problems and $3{,}000$ for 1D, generated using Hammersley sampling~\cite{hammersley1960monte}.

\subsection{Convergence properties of multi-preconditioned LBFGS}
First, we investigate the impact of the scaling strategies discussed in Sect.~\ref{sec:multi_prec} 
on the performance of the MP-LBFGS method.
Fig.~\ref{fig:plotscalings_robustness} 
shows that, independently of the chosen decomposition, the UniS strategy leads to unstable training.
We further observe that as the number of subdomains increases (e.g., from $2\times2$ to $4\times4$), the SPM strategy significantly outperforms the LSS approach. 
This behavior can be attributed to the increasing distance between successive iterates after each line-search step.
We also observe that, even for a fixed number of subdomains, the choice of decomposition strongly affects MP-LBFGS performance.
For example, moving from a $2\times2$ to a $4\times1$ decomposition significantly degrades the LSS method, whereas the SPM strategy achieves an error reduction of more than one order of magnitude.

\begin{figure}[t]
    \centering
    \includegraphics{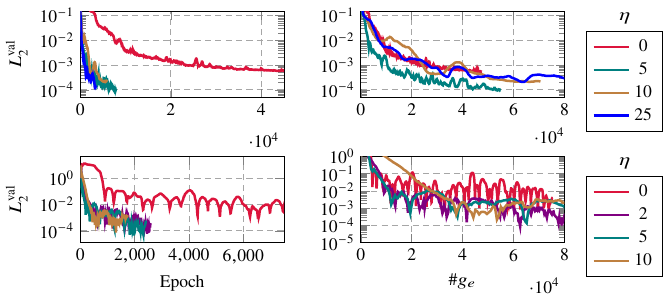}

    \caption{Convergence of MP-LBFGS method with respect to varying number of local LBFGS iterations ($\eta$) for Burgers (top, 8 subdomains) and 1D Poisson (bottom, 20 subdomains), where the value $\eta=0$ corresponds to standard LBFGS method.}
    \label{fig:expII2x4}
\end{figure}

Next, we compare the performance of MP-LBFGS  with standard LBFGS.
Results are reported in terms of epochs and effective parallel gradient evaluations $(\# g_e)$.
The metric $\# g_e$ approximately reflects the per-device computational work, while excluding communication overhead, costs of parameter updates and loss evaluations (see 
Tab.~\ref{tab:computational_cost}). 
The number of epochs serves as a proxy for communication cost.
\newadd{Here, we also note that, for all presented experiments, the number of Newton iterations in the SPM method, denoted as $\#\text{its}_{\text{Newton}}$, is on average equal to two.}

We first examine the performance of MP-LBFGS while varying the number of local LBFGS iterations $\eta$, with $\eta=0$ corresponding to standard LBFGS.
As shown in Fig.~\ref{fig:expII2x4},  
all MP-LBFGS configurations converge significantly faster in terms of epochs, thereby reducing communication cost.
Moreover, MP-LBFGS requires a comparable, and in several cases lower, number of effective gradient evaluations $(\# g_e)$.
We also highlight, that for several values of $\eta$, MP-LBFGS achieves up to an order-of-magnitude reduction in the validation error $L_2^{\mathrm{val}}$.

\begin{figure}[t]
  \centering
  \includegraphics{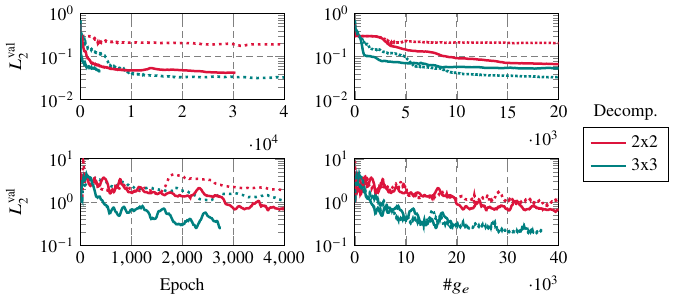}

  \caption{
    Convergence of MP-LBFGS (solid lines, $\eta=5$) and LBFGS (dotted lines) for Burgers (top) and Poisson (bottom) problem.
    \newadd{Experiments terminate after $\# g_e$ of $20,000$ or $40,000$ is reached.}}
  \label{fig:expIIIB}
\end{figure}

Second, we compare LBFGS and MP-LBFGS with $\eta=5$ on the 2D Poisson and Burgers problems using $2\times2$ and $3\times3$ decompositions, while keeping the number of collocation points per subdomain fixed.
As shown in Fig.~\ref{fig:expIIIB}, 
MP-LBFGS consistently reduces the number of epochs while requiring comparable $\# g_e$ and achieving comparable accuracy across different decompositions.
\newadd{Note that, here, one epoch of the MP-LBFGS algorithm with $\eta = 5$ corresponds to seven gradient evaluations; see Tab.~\ref{tab:computational_cost} 
for details.}    
The only exception is the Burgers problem with a $3\times3$ decomposition, where LBFGS attains higher accuracy. \newdel{; this outlier will be investigated in future work.}
\newadd{Our empirical evidence suggests that, for this particular example, the non-convexity of the problem causes MP-LBFGS to converge to a local minimum with a larger error than L-BFGS.}

\section{Conclusion}
\label{sec:conclusion}
In this work, we developed a multi-preconditioned LBFGS (MP-LBFGS) framework for training FBPINNs.
The proposed approach exploits the FBPINN architecture to compute local LBFGS search directions in parallel, while a global right-preconditioned LBFGS iteration ensures consistent convergence of the global minimization problem.
A key component of the method is a novel subspace minimization strategy (SPM) for scaling locally computed search directions.
Numerical experiments suggest that the proposed MP-LBFGS can achieve faster convergence, and higher accuracy than the standard LBFGS, while maintaining lower communication overhead.
Future work will focus on parallel implementation, \newadd{allowing for wall-clock comparison and the incorporation of more complex examples as well as on} incorporation of coarse spaces\newdel{,}\newadd{;} and on reducing the cost of solving the SPM subproblem.

\vspace{-0.3cm}
\begin{acknowledgement}
The work of A.Ks, S.G. and A.Ko. benefited from ANITI, funded by the  France 2030 program under Grant Agreement No.~ANR-23-IACL-0002. 
The numerical results were carried out using HPC resources from GENCI-IDRIS (Grant No.~AD011015766R1).
\end{acknowledgement}

\vspace{-1cm}
\bibliographystyle{spmpsci}
\bibliography{Biblio}

\end{document}